\documentclass[a4paper]{article}
\usepackage{amsmath}
\usepackage[psamsfonts]{amssymb}

\newcommand{\cC}{{\cal C}}
\newcommand{\Ci}{C^{\infty}}
\newcommand{\cD}{{\cal D}}
\newcommand{\cL}{{\cal L}}

\newcommand{\p}{\partial}

\newcommand{\cS}{{\cal S}}

\newcommand{\cA}{{\cal A}}
\newcommand{\cX}{{\cal X}}

\newcommand{\cM}{{\cal M}}
\newcommand{\cI}{{\cal I}}
\newcommand{\cZ}{{\cal Z}}
\newcommand{\raa}{\rightarrow}
\newcommand{\E}{\ell}
\newcommand{\m}{\!\!\mid}
\newcommand{\N}{{\mathbf N}}

\newcommand{\R}{{\mathbf R}}

\newcommand{\lp}{\left(}
\newcommand{\rp}{\right)}

\newcommand{\ap}{\alpha}

\newcommand{\La}{\lambda}

\newcommand{\DE}{\Delta}
\newcommand{\si}{\sigma}
\newcommand{\z}{\zeta}
\newcommand{\be}{\begin{equation}}
\newcommand{\ee}{\end{equation}}

\newtheorem{theo}{Theorem}
\newtheorem{lem}{Lemma}

\begin{document}
\title{Lie algebraic characterization of
manifolds\thanks{The research of Janusz Grabowski was supported by
the Polish Ministry of Scientific Research and Information
Technology under the grant No. 2 P03A 020 24, that of Norbert
Poncin by grant C.U.L./02/010.}}\author{Janusz
Grabowski\footnote{E-mail: jagrab@impan.gov.pl}\,, Norbert
Poncin\footnote{E-mail: norbert.poncin@uni.lu}} \maketitle
\begin{abstract}
Results on characterization of manifolds in terms of certain Lie
algebras growing on them, especially Lie algebras of differential
operators, are reviewed and extended. In particular, we prove that
a smooth (real-analytic, Stein) manifold is characterized by the
corresponding Lie algebra of linear differential operators, i.e.
isomorphisms of such Lie algebras are induced by the appropriate
class of diffeomorphisms of the underlying manifolds.

{\bf Keywords}: algebraic characterization; smooth, real-analytic,
and Stein manifolds; automorphisms; Lie algebras; differential
operators; principal symbols

{\bf MSC (2000)}: 17B63 (primary), 13N10, 16S32, 17B40, 17B65,
53D17 (secondary)
\end{abstract}

\section{Algebraic characterizations of manifold
structures}

Algebraic characterizations of topological spaces and manifolds
can be traced back to the work of I.~Gel'fand and A.~Kolmogoroff
\cite{GK} in which compact topological spaces $K$ are
characterized by the algebras $\cA=C(K)$ of continuous functions
on them. In particular, points $p$ of these spaces are identified
with maximal ideals $p^*$ in these algebras consisting of
functions vanishing at $p$. This identification easily implies
that isomorphisms of the algebras $C(K_1)$ and $C(K_2)$ are
induced by homeomorphisms between $K_1$ and $K_2$, so that the
algebraic structure of $C(K)$ characterizes $K$ uniquely up to
homeomorphism. Here $C(K)$ may consist of complex or real
functions as well.

All above can be carried over when we replace $K$ with a compact
smooth manifold $M$ and $C(K)$ with the algebra $\cA=C^\infty(M)$
of all real smooth functions on $M$. In this case algebraic
isomorphisms between $C^\infty(M_1)$ and $C^\infty(M_2)$ are
induced by diffeomorphisms $\phi:M_2\rightarrow M_1$, i.e. they
are of the form $\phi^*(f)=f\circ\phi$. If our manifold is
non-compact then no longer maximal ideals must be of the form
$p^*$. However, there is still an algebraic characterization of
ideals $p^*$ as one-codimensional (or just maximal
finite-codimensional) ones, so that the algebraic structure of
$C^\infty(M)$ still characterizes $M$, for $M$ being Hausdorff and
second countable. A similar result is true in the real-analytic
(respectively the holomorphic) case, i.e. when we assume that $M$
is a real-analytic (respectively a Stein) manifold and that $\cA$
is the algebra of all real-analytic (respectively all holomorphic)
functions on $M$ (see \cite{JG}). As it was recently pointed out
to us by Alan Weinstein, the assumption that $M$ is second
countable is crucial for the standard proofs that
one-codimensional ideals are of the form $p^*$. This remark
resulted in alternative proofs \cite{Mrc, Gra03}, which are valid
without additional assumptions about manifolds.

There are other algebraic structures canonically associated with a
smooth manifold $M$, for example the Lie algebra $\cX(M)$ of all
smooth vector fields on $M$ or the associative (or Lie) algebra
$\cD(M)$ of linear differential operators acting on $C^\infty(M)$.
Characterization of manifolds by associated Lie algebras is a
topic initiated in $1954$ with the well-known paper \cite{PS} of
L.E. Pursell and M.E. Shanks that appeared in the Proceedings of
the American Mathematical Society.

The main theorem of this article states that, if the Lie algebras
$\cX_c(M_i)$ of smooth compactly supported vector fields of two
smooth manifolds $M_i$ ($i\in\{1,2\}$) are isomorphic Lie
algebras, then the underlying manifolds are diffeomorphic,
and---of course---vice versa. The central idea of the proof of
this result is the following. If $M$ is a smooth variety and $p$ a
point of $M$, denote by $\cI(\cX_c(M))$ the set of maximal ideals
of $\cX_c(M)$ and by $p^\infty$ the maximal ideal
\[p^\infty=\{X\in\cX_c(M):X\mbox{ is flat at }p\}.\] The map $p\in M\raa
p^\infty\in\cI(\cX_c(M))$ being a bijection and the property "maximal
ideal" being an algebra-isomorphism invariant, the correspondence
\[\begin{array}{ccc}p_1^\infty\in\cI(\cX_c(M_1))&\raa&p_2^\infty\in\cI(\cX_c(M_2))\\
\uparrow&&\downarrow\\p_1\in
M_1 & & p_2\in M_2\end{array}\] is a bijection and even a
diffeomorphism.

Similar upshots exist in the analytic cases. Note that here "flat"
means zero and that the Lie algebra of all $\R$-analytic vector
fields of an $\R$-analytic compact connected manifold is simple,
so has no proper ideals. Hence, the above maximal ideals in
particular, and ideals in general, are of no use in these cases.
Maximal finite-codimensional subalgebras turned out to be an
efficient substitute for maximal ideals. This idea appeared around
$1975$ in several papers, e.g. in \cite{Am}, and is the basis of
the general algebraic framework developed in \cite{JG} and
containing the smooth, the $\R$-analytic and the holomorphic
cases. Here, if $\cM(\cX_{\bullet}(M))$ is the set of all maximal
finite-codimensional subalgebras of $\cX_{\bullet}(M)$, where
subscript $\bullet$ means smooth, $\R$-analytic or holomorphic,
the fundamental bijection is $p\in M\raa
p^0\in\cM(\cX_{\bullet}(M))$ with
\[p^0=\{X\in\cX_{\bullet}(M):(Xf)(p)=0,\forall f\in C^{\bullet}(M)\}.\]

This method works well for a large class of the Lie algebras of
vector fields which are simultaneously modules over the algebra of
functions $\cA$. A pure algebraic framework in this direction has
been developed by S.~M.~Skryabin \cite{S} who has proven a very
general "algebraic Pursell-Shanks theorem":

\begin{theo}If $\cL_i$ is a Lie subalgebra in the Lie algebra $Der\cA_i$ of
derivations of a commutative associative unital algebra $\cA_i$ over a
field of characteristic $0$ which is simultaneously an $\cA_i$-submodule
of $Der\cA_i$ and is non-singular in the sense that $\cL_i(\cA_i)=\cA_i$,
$i=1,2$, then any Lie algebra isomorphism $\Phi:\cL_1\rightarrow\cL_2$ is
induced by an associative algebra isomorphism
$\phi^*:\cA_2\rightarrow\cA_1$, i.e.
$$\Phi(X)(f)=(\phi^*)^{-1}(X(\phi^*(f))).$$
\end{theo}
Note that Skryabin's proof does not refer to the structure of
maximal ideals in $\cA$ but uses the $\cA$-module structure on
$\cL$.

Other types of Lie algebras of vector fields (see e.g. \cite{KMO})
have also been considered but the corresponding methods have been
developed for each case separately. Let us mention the Lie
algebras of vector fields preserving a given submanifold
\cite{Ko}, a given (generalized) foliation \cite{JG3}, a
symplectic or contact form \cite{O}, the Lie algebras of
Hamiltonian vector fields or Poisson brackets of functions on a
symplectic manifold \cite{AG} and Jacobi brackets in general
\cite{JG2}, Lie algebras of vector fields on orbit spaces and
$G$-manifolds \cite{A}, Lie algebras of vector fields on affine
and toric varieties \cite{HM,CGM, Si}, Lie algebroids \cite{GG},
and many others.

In our work \cite{GP}, we have examined the Lie algebra $\cD(M)$
of all linear differential operators on the space $\Ci(M)$ of
smooth functions of $M$, its Lie subalgebra $\cD^1(M)$ of all
first-order differential operators and the Poisson algebra
$\cS(M)=Pol(T^*M)$ of all polynomial functions on the cotangent
bundle, the symbols of the operators in $\cD(M)$. We have obtained
in each case Pursell-Shanks type results in a purely algebraic
way. Furthermore, we have provided an explicit description of all
the automorphisms of any of these Lie algebras.

In this notes we depict this last paper assuming a philosophical
and pedagogic point of view, we prove a general algebraic
Pursell-Shanks type result, and, with the help of some topology,
we extend our smooth Pursell-Shanks type result from \cite{GP} to
the real-analytic and holomorphic cases.

\section{Lie algebras of differential operators}

\subsection{Abstract definitions}

The goal being a work on the algebraic level, we must define some
general algebra, modelled on $\cD(M)$, call it a quantum Poisson
algebra. A quantum Poisson algebra is an associative
\textit{filtered} algebra $\cD=\cup_i\cD^{i}$ ($\cD^{i}=\{0\},$
for $i<0$) with unit $1$,
\[\cD^{i}\cdot\cD^{j}\subset\cD^{i+j},\] such that the canonical
\textit{Lie bracket} verifies
\begin{equation}[\cD^{i},\cD^j]\subset\cD^{i+j-1}.\label{Liefil}
\end{equation}
Note that ${\cal A}=\cD^0$ is an associative commutative
subalgebra of $\cD$ and that $K$, the underlying field, is
naturally imbedded in $\cA$.

Similarly we heave the algebra $\cS(M)=Pol(T^*M)$ (or $\cS(M)=\Gamma(\cS
TM)$) of smooth functions on $T^*M$ that are polynomial along the fibers
(respectively, of symmetric contravariant tensor fields of $M$), classical
counterpart of $\cD(M)$, on the algebraic level and define a classical
Poisson algebra as a \textit{commutative} associative \textit{graded}
algebra $\cS=\oplus_i\cS_i$ ($\cS_i=\{0\},$ for $i<0$) with unit $1$,
\[\cS_i\cS_j\subset \cS_{i+j},\] which is equipped with a
\textit{Poisson} bracket $\{.,.\}$ such that
\[\{\cS_i,\cS_j\}\subset \cS_{i+j-1}.\]
Here ${\cal A}=\cS_0$ is obviously an associative and
Lie-commutative subalgebra of $\cS$.

Let us point out that quantum Poisson algebras canonically induce
classical Poisson algebras. Indeed, starting from
$\cD=\cup_i\cD^{i}$, we get a graded vector space when setting
$\cS_i=\cD^{i}/\cD^{i-1}$. If the degree $deg(D)$ of an arbitrary
non-zero $D\in\cD$ is given by the lowest filter that contains $D$
and if $cl_i$ denotes the class in the quotient $\cS_i$, we define
the principal symbol $\sigma(D)$ of $D$ by
\[\sigma(D)=\mbox{ cl}_{deg(D)}(D)\]
and the symbol $\sigma_i(D)$ of order $i\ge deg(D)$ by
\[\sigma_i(D)=\mbox{
cl}_i(D)=\begin{cases} 0,\mbox{ if }i>deg(D), &\\
\sigma(D),\mbox{ if }i=deg(D).&\end{cases}\] Now the
\textit{commutative} associative multiplication and the
\textit{Poisson} bracket are obtained as follows:
\[\si(D_1)\si(D_2)=\si_{deg(D_1)+deg(D_2)}(D_1\cdot D_2)\;\;(D_1,D_2\in\cD)\]and
\[\{\si(D_1),\si(D_2)\}=\si_{deg(D_1)+deg(D_2)-1}([D_1,D_2])\;\;(D_1,D_2\in\cD).\]
Remark that the commutativity of the associative multiplication is
a direct consequence of Equation (\ref{Liefil}) and that Leibniz's
rule simply passes from $[.,.]$ to $\{.,.\}$. So "dequantization"
is actually a passage from non-commutativity to commutativity, the
trace of non-commutativity on the classical level being the
Poisson bracket. Note also that the quantum Poisson algebra $\cD$
and the induced classical limit $\cS$ of $\cD$ have the same basic
algebra $\cA$. The principal symbol map $\sigma:\cD\rightarrow\cS$
has the following important property:
\[
\{\sigma(D_1),\sigma(D_2)\}=\begin{cases} \sigma([D_1,D_2])&\\
\mbox{or}&\\
0&\end{cases}\ .\]

There is a canonical quantum Poisson algebra associated with any unital
associative commutative algebra $\cA$, namely the algebra $\cD(\cA)$ of
linear differential operators on $\cA$. Note that this algebraic approach
to differential operators goes back to some ideas of Grothendieck and that
it was extensively developed by A.~M.~Vinogradov. The filtration
$\cD^i(\cA)$ is defined inductively: $\cD^0(\cA)=\{ m_f:f\in\cA\}$, where
$m_f(g)=f\cdot g$ (so $\cD^0(\cA)$ is canonically isomorphic with $\cA$),
and
$$\cD^{i+1}(\cA)=\{ D\in Hom(\cA):[D,m_f]\in\cD^{i}(\cA), \textit{ for all
} f\in\cA\},$$ where $[.,.]$ is the commutator.

It can be seen that, in the fundamental example $\cA=\Ci(M)$, we
get $\cD=\cD(\cA)=\cD(M), \cS=\cS(M)=Pol(T^*M)=\Gamma(\cS TM)$ and
the above algebraically defined Poisson bracket coincides with the
canonical Poisson bracket on $\cS$, i.e. the standard symplectic
bracket in the first interpretation of $\cS$, and the symmetric
Schouten bracket in the second. The situation is completely
analogous in real-analytic and holomorphic cases. Note only that
in the holomorphic case the role of $T^*M$ is played by a complex
vector bundle $T^*_{(1,0)}M$ over $M$ whose sections are
holomorphic 1-forms of the type $(1,0)$ and the Poisson structure
on $T^*_{(1,0)}M$ is represented by a holomorphic bivector field
being a combination of wedge products of vector fields of type
$(1,0)$.

\subsection{Key-idea}

Remember that the first objective is to establish Pursell-Shanks
type results, that is---roughly speaking---to deduce a geometric
conclusion from algebraic information. As in previous papers on
this topic, functions should play a central role. So our initial
concern is to obtain an algebraic characterization of "functions",
i.e., in the general algebraic context of an arbitrary quantum
Poisson algebra $\cD$, of ${\cal A}\subset\cD$, and more generally
of all filters $\cD^{i}\subset\cD$.

It is clear that \[\begin{array}{lll}{\cal
A}&\subset&\mbox{Nil}(\cD)\\&:=&\{D\in\cD:\forall\DE\in\cD,\exists
n\in\N:\stackrel{n}{\overbrace{[D,[D,\ldots[D}},\DE]]]=0\}\end{array}\]
and that \[\cD^{i+1}\subset\{D\in\cD:[D,{\cal A}]\subset \cD^{i}\}\mbox{
}(i\ge -1).\] Our conjecture is that functions (respectively, $(i+1)$th
order "differential operators") are the only locally nilpotent operators,
i.e. the sole operators whose repeated adjoint action upon any operator
ends up by zero (respectively, the only operators for which the commutator
with all functions is of order $i$).

It turns out that both guesses are confirmed if we show that for
any $D\in\cD$,
\[\forall f\in{\cal A},\exists
n\in\N:\stackrel{n}{\overbrace{[D,[D,\ldots[D}},f]]]=0\Longrightarrow
D\in{\cal A}.\] We call this property, which states that if the
repeated adjoint action of an operator upon any function ends up
by zero then this operator is a function, the distinguishing
property of the Lie bracket.

At this point it is natural to ask if any bracket is
distinguishing and---in the negative---if the commutator bracket
of $\cD(M)$ is. Obviously, the algebra of all linear differential
operators acting on the polynomials in a variable $x\in\R$ is a
simple example of a non-distinguishing Lie algebra. It suffices to
consider the operator $d/dx$. The answer to the second question is
positive. We refrain from describing here the technical
constructive proof given in \cite{GP}. We will present further
another topological proof, which can be adapted to real-analytic
and holomorphic cases.

\subsection{An algebraic Pursell-Shanks type result}

We aspire to give an algebraic proof of the theorem stating that
two manifolds $M_1$ and $M_2$ are diffeomorphic if the Lie
algebras $\cD(M_1)$ and $\cD(M_2)$ are isomorphic. So let us
consider a Lie algebra isomorphism
\[\Phi:\cD_1\raa\cD_2\]
between two quantum Poisson algebras $\cD_1$ and $\cD_2$.

In the following, we discuss two necessary assumptions.

\subsubsection{Distinguishing property}

The next proposition is a first step towards our aim.

\medskip
If $\cD_1,\cD_2$ are \textit{distinguishing} quantum Poisson algebras then
$\Phi$ respects the filtration.

\medskip\noindent
The proof is by induction on the "order of differentiation" and
uses the above algebraic characterizations of functions and
filters, hence the distinguishing character of $\cD_1$ and
$\cD_2$. For instance, $\Phi({\cal
A}_1)=\Phi(\mbox{Nil}(\cD_1))=\mbox{Nil}(\cD_2)={\cal A}_2$ (so
that in particular $\Phi({\cal D}^0_1)\subset{\cal D}^0_2$).

\subsubsection{Non-singularity property}

Since $\Phi({\cal A}_1)={\cal A}_2$, the Lie algebra isomorphism
$\Phi$ restricts to a vector space isomorphism between ${\cal
A}_1$ and ${\cal A}_2$. If this restriction respected the
associative multiplication it would be an associative algebra
isomorphism, which---as well known---would in the geometric
context, $\cD_i=\cD(M_i)$ ($i\in\{1,2\}$), be induced by a
diffeomorphism between $M_1$ and $M_2$.

So the question arises if we are able to deduce from the Lie
algebra structure any information regarding the associative
algebra structure and in particular the left and right
multiplications $\E_f:\cD\ni D\raa f\cdot D\in\cD$ and $r_f:\cD\ni
D\raa D\cdot f\in\cD$ by a function $f\in\cA$.

Observe first that $\E_f$ and $r_f$ commute with the adjoint
action by functions, i.e. are members of the centralizer of
$ad\,\cA$ in the Lie algebra $End(\cD)$ of endomorphisms of $\cD$,
which is the Lie subalgebra $\cC(\cD)=\{\Psi\in End(\cD):\Psi\circ
ad\,{\cal A} =ad\,{\cal A}\circ\Psi\}$.

On the other hand, it is not possible to extract from the Lie bracket more
information than
\begin{equation}\E_f-r_f=ad\,f,\label{central}\end{equation}
where the right hand side is of course a \textit{lowering} member
of the centralizer, i.e. a mapping in the centralizer which lowers
the order of differential operators.

Thus the centralizer might be the brain wave. In particular it
should, in view of (\ref{central}), be possible to describe it as
the algebra of those endomorphisms $\Psi$ of $\cD$ that respect
the filtration and are of the form $\Psi = \E_f+\Psi_1,$ where
$f=\Psi(1)$ and $\Psi_1\in\cC(\cD)$ is lowering. When trying to
prove this conjecture, we realize that it holds if
$[\cD^1,\cA]=\cA$, in the sense that any function is a finite sum
of brackets. In the geometric context this means that
$[\cX(M),\Ci(M)]=\cX(M)\lp\Ci(M)\rp=\Ci(M)$, a non-singularity
assumption that appears in many papers of this type, e.g.
\cite{JG,S}, and is of course verified. Hence, a second
proposition:

If $\cD$ is a \textit{non-singular} and distinguishing quantum
Poisson algebra then any $\Psi\in\cC(\cD)$ respects the filtration
and $\Psi=\E_{\Psi(1)}+\Psi_1,$ $\Psi_1\in\cC(\cD)$ being
lowering.

\subsubsection{Isomorphisms}

Having in view to use the centralizer to show that $\Phi$ respects the
associative multiplication, we must visibly read the Lie algebra
isomorphism $\Phi:\cD_1\raa\cD_2$ as Lie algebra isomorphism
$\Phi_*:End(\cD_1)\raa End(\cD_2)$. We only need set
$\Phi_*(\Psi)=\Phi\circ\Psi\circ\Phi^{-1}$, $\Psi\in End(\cD_1)$. As
$\Phi_*(\cC(\cD_1))=\cC(\cD_2)$, it follows from the above depicted
structure of the centralizer that for any $f,g\in{\cal A}_1$,
$\Phi_*(\E_f)(\Phi(g))=\lp\Phi_*(\E_f)\rp(1)\cdot\Phi(g),$ i.e.
\be\label{z}\Phi(f\cdot g)=\Phi(f\cdot\Phi^{-1}(1))\cdot\Phi(g).\ee

Remark that $\z=\Phi^{-1}(1)\in\cZ(\cD_1)$, where $\cZ(\cD_1)$
denotes the center of the Lie algebra $\cD_1$. In view of
(\ref{z}), \be\label{z1}\Phi(f\cdot g)=\Phi(g\cdot
f)=\Phi(f)\cdot\Phi(g\cdot\z)=\Phi(f)\cdot\Phi(\zeta\cdot
g)=\Phi(\z^2)\cdot \Phi(f)\cdot\Phi(g), \ee for all $f,g\in\cA_1$.
This in turn implies that $\Phi(\z^2)$ is invertible, say
$\Phi(\z^2)=\kappa^{-1}$, since the image of $\Phi$, so $\cA_2$,
is contained in the ideal generated by $\Phi(\z^2)$. If we put
$A(f)=\kappa^{-1}\cdot\Phi(f)$, then, due to (\ref{z1}), $A(f\cdot
g)=A(f)\cdot A(g)$, so $A$ is an associative algebra isomorphism.
Thus we get the following algebraic Pursell-Shanks type theorem.

\begin{theo}\label{APS} Let $\cD_i$ $($$i\in\{1,2\}$$)$ be
non-singular and distinguishing quantum Poisson algebras. Then every Lie
algebra isomorphism $\Phi:\cD_1\raa\cD_2$ respects the filtration and its
restriction $\Phi\m_{{\cal A}_1}$ to ${\cal A}_1$ has the form
$\Phi\m_{{\cal A}_1}=\kappa A,$ where $\kappa\in\cA_2$ is invertible and
central in $\cD_2$ and $A:{\cal A}_1\raa {\cal A}_2$ is an associative
algebra isomorphism.\label{PS}
\end{theo}

{\bf Remark.} The central elements in $\cD(M)$ are just constants.
This immediately follows from the well-known corresponding
property of the symplectic Poisson bracket on
$\cS=\cS(M)=Pol(T^*M)$. If this symplectic property holds good for
the classical limit $\cS$ of $\cD$, we say that $\cD$ is
\textit{symplectic}. We have assumed this property to obtain a
version of Theorem \ref{APS} in \cite{GP}. Now we have proven that
this assumption is superfluous.

\subsection{Isomorphisms of the Lie algebras of differential operators}

Let now $M$ be a finite-dimensional paracompact and second
countable smooth (respectively real-analytic, Stein) manifold, let
$\cA=\cA(M)$ be the commutative associative algebra of all real
smooth (respectively real-analytic, holomorphic) functions on $M$,
let $\cD=\cD(M)=\cD(\cA(M))$ be the corresponding algebra of
differential operators, and let $\cS=\cS(M)$ be the classical
limit of $\cD(M)$.

\begin{lem} The quantum Poisson algebra $\cD(M)$ is distinguishing and
non-\linebreak singular.
\end{lem}

{\it Proof.} Let us work first in the smooth case. Let $D$ be a
linear differential operator on $\cA(M)$ such that for every
$f\in\cA(M)$ there is $n$ for which $(ad_D)^n(f)=0$, where
$ad_D(D')=[D,D']$ is the commutator in the algebra of differential
operators. It suffices to show that $D\in\cA(M)$.

The algebra $\cA(M)$ admits a complete metric, which makes it into
a topological algebra such that all linear differential operators
are continuous (see section 6 of the book \cite{KM}). It is then
easy to see that $Ker((ad_D)^n)=\{ f\in\cA(M):(ad_D)^n(f)=0\}$,
$n=1,2,\dots$, is a closed subspace of $\cA(M)$. By assumption,
$\bigcup_nKer((ad_D)^n)=\cA(M)$, so $Ker((ad_D)^{n_0})=\cA(M)$ for
a certain $n_0$ according to the Baire property of the topology on
$\cA(M)$. Passing now to principal symbols, we can write
$X_D^{n_0}(f)=0$ for all $f\in \cA(M)$, where $X_D$ is the
Hamiltonian vector field on $T^*M$ of the principal symbol
$\sigma(D)$ of $D$ with respect to the canonical Poisson bracket
on $T^*M$. Here we regard $\cA(M)$ as canonically embedded in the
algebra of polynomial functions on $T^*M$. Hence the 0-order
operator $(ad_{g})^{n_0}(X_D^{n_0})$, which is the multiplication
by $n_0!{(X_D(g))^{n_0}}$, vanishes on $\cA(M)$ for all
$g\in\cA(M)$, so $X_D(g)=\{\sigma(D),g\}=0$ for all $g\in\cA(M)$.
This in turn implies that $\sigma(D)\in\cA(M)$, so $D\in\cA(M)$.

\medskip
The proof in the holomorphic case is completely analogous (for the
topology we refer to section 8 of \cite{KM}), so let us pass to
the real-analytic case. Now, the natural topology is not
completely metrizable and we cannot apply the above procedure.
However, there is a Stein neighbourhood $\tilde M$ of $M$, i.e. a
Stein manifold $\tilde M$ of complex dimension equal to the
dimension of $M$, containing $M$ as a real-analytic closed
submanifold, whose germ along $M$ is unique \cite{Grau}, so that
$D$ can be complexified to a linear holomorphic differential
operator $\tilde D$. The complexified operator has clearly the
analogous property: for every holomorphic function $f$ on $\tilde
M$ there is $n$ for which $(ad_{\tilde D})^n(f)=0$, so $\tilde D$,
thus $D$, is a multiplication by a function.

The non-singularity of $\cD(M)$ follows directly from Proposition 3.5 of
\cite{JG}. \rule{1.5mm}{2.5mm}

\medskip
The above Lemma shows that we can apply Theorem \ref{APS} to $\cD(M)$ in
all, i.e. smooth, real-analytic, and holomorphic, cases. We obtain in the
same way---mutatis mutandis---Pursell-Shanks type results for the Lie
algebras $\cS(M)$ and $\cD^1(M)$:

\begin{theo}The Lie algebras $\cD(M_1)$ and $\cD(M_2)$ $($respectively $\cS(M_1)$
and $\cS(M_2)$, or $\cD^1(M_1)$ and $\cD^1(M_2)$$)$ of all
differential operators $($respectively all  symmetric
contravariant tensors, or all differential  operators  of order
$1$$)$ on two smooth (respectively real-analytic, holomorphic)
manifolds $M_1$ and $M_2$ are isomorphic if and only if the
manifolds $M_1$ and $M_2$ are smoothly (respectively
bianalytically, biholomorphically) diffeomorphic.
\end{theo}

\subsection{Automorphisms of the Lie algebras of differential operators}

In view of the above Pursell-Shanks type results, the study of the
Lie algebra isomorphisms, e.g. between $\cD(M_1)$ and $\cD(M_2)$,
can be reduced to the examination of the corresponding Lie algebra
automorphisms, say $\Phi\in Aut(\cD(M),[.,.])$. The standard idea
in this kind of problems is the simplification of the considered
arbitrary automorphism $\Phi$ via multiplication by automorphisms
identified a priori. Here, the automorphism $A=\kappa^{-1}\Phi\in
Aut(\Ci(M),\cdot)$ induced by $\Phi$ (see Theorem \ref{PS}) can
canonically be extended to an automorphism $A_*\in
Aut(\cD(M),[.,.])$. It suffices to set $A_*(D)=A\circ D\circ
A^{-1}$, for any $D\in\cD(M)$. Then, evidently,
\[\Phi_1=A_*^{-1}\circ\Phi\in Aut(\cD(M),[.,.])\mbox{  and  }\Phi_1\m_{\Ci(M)}=\kappa\,id,\]
where $id$ is the identity map of $\Ci(M)$.

Some notations are necessary. In the following, we use the
canonical splitting $\cD(M)=\Ci(M)\oplus\cD_c(M)$, where
$\cD_c(M)$ is the Lie algebra of differential operators vanishing
on constants. Moreover, we denote by $\bullet_{0}$ and
$\bullet_{c}$ the projections onto $\Ci(M)=\cD^0(M)$ and
$\cD_c(M)$ respectively.

\subsubsection{Automorphisms of $\mathbf{\cD^1(M)}$}

The formerly explained reduction to the problem of the
determination of all automorphisms $\Phi_1,$ which coincide with
$\kappa\,id$ on functions, is still valid. Furthermore, in the
case of the Lie algebra $\cD^1(M)$ the preceding splitting reads
$\cD^1(M)=\Ci(M)\oplus\cD_c^1(M)=\Ci(M)\oplus\cX(M)$.

It follows from the automorphism property that
\[(\Phi_1)_c\m_{\cX(M)}=id\mbox{  and  }(\Phi_1)_0\m_{\cX(M)}\in {\cal
Z}^1(\cX(M),\Ci(M)),\] where ${\cal Z}^1(\cX(M),\Ci(M))$ is the
space of $1$-cocycles of the Lie algebra of vector fields
canonically represented upon functions by Lie derivatives. We know
that these cocycles are locally given by
\begin{equation}(\Phi_1)_0\m_{\cX(M)}=\La\,div+df,\label{cocycles}\end{equation} where $\La\in\R$ and $f$
is a smooth function. When trying to globally define the right
hand side of Equation (\ref{cocycles}), we naturally substitute a
closed $1$-form $\omega\in\Omega^1(M)\cap ker\,d$ to the exact
$1$-form $df$. In order to globalize the divergence $div$, note
the following. If $M$ is oriented by a volume $\Omega$, we have
$div_{\Omega}=div_{-\Omega}$. So the divergence can be defined
with respect to $\mid\Omega\mid$. This pseudo-volume may be viewed
as a pair $\{\Omega,-\Omega\}$ and exists on any manifold,
orientable or not. Alternatively, we may interpret
$\mid\Omega\mid$ as a nowhere vanishing tensor density of weight
$1$, i.e. as a section of the vector bundle $I\!\!F_1(TM)$ of
$1$-densities, which is everywhere non-zero. This bundle being of
rank $1$, the existence of such a section is equivalent to the
triviality of the bundle. However, the fiber bundle $I\!\!F_1(TM)$
is known to be trivial for any manifold $M$. Thus, a nowhere
vanishing tensor $1$-density $\rho_0$ always exists and the
divergence can be defined with respect to this $\rho_0$. For a
more rigorous approach the reader is referred to \cite{GP}.
Eventually, the above cocycles globally read
\[(\Phi_1)_0\m_{\cX(M)}=\La\,div_{\rho_0}+\omega.\]

Let us fix a divergence $div$. For any
$f+X\in\cD^1(M)=\Ci(M)\oplus\cX(M)$, we then obtain
\[\Phi_1(f+X)=\kappa
f+\La\,div\,X+\omega(X)+X.\] Since the initial arbitrary
automorphism $\Phi\in Aut(\cD^1(M),[.,.])$ has been decomposed as
$\Phi=A_*\circ\Phi_1$, with $A_*$ induced by a diffeomorphism
$\varphi\in Diff(M)$ and denoted $\varphi_*$ below, we finally
have the theorem:

\begin{theo} A linear map $\Phi:\cD^1(M)\raa\cD^1(M)$ is an automorphism
of the Lie algebra $\cD^1(M)=\Ci(M)\oplus\cX(M)$ of linear
first-order differential operators on $\Ci(M)$ if and only if it
can be written in the form
\[\Phi(f+X)=\lp\kappa
f+\La\,div\,X+\omega(X)\rp\circ\varphi^{-1}+\varphi_*(X),\] where
$\varphi$ is a diffeomorphism of $M$, $\lambda$, $\kappa$ are
constants $($$\kappa\ne 0$$)$, $\omega$ is a closed $1$-form on
$M$, and $\varphi_*$ is defined by
\[(\varphi_*(X))(f)=(X(f\circ\varphi))\circ\varphi^{-1}.\] All the objects $\varphi,
\lambda,\kappa,\omega$ are uniquely determined by $\Phi$.
\label{d1}
\end{theo}
\subsubsection{Automorphisms of $\mathbf{\cD(M)}$}

As previously, we need only seek the automorphisms $\Phi_1\in
Aut(\cD(M),[.,.])$, such that $\Phi_1(f)=\kappa f$, $f\in\Ci(M)$.
Such an automorphism visibly restricts to a similar automorphism
of the Lie algebra $\cD^1(M).$ Hence,
\begin{equation}\Phi_1(f+X)=\kappa
f+\La\,div\,X+\omega(X)+X.\label{d1f1}\end{equation}

Since $\omega\in\cZ^1(\cX(M),\Ci(M))$, it is reasonable to think
that $\omega$ might be extended to
$\overline{\omega}\in\cZ^1(\cD(M),\cD(M))=Der(\cD(M))$, where
$Der(\cD(M))$ is the Lie algebra of all derivations of
$(\cD(M),[.,.])$. If in addition this derivation
$\overline{\omega}$ were lowering, it would generate an
automorphism $e^{\overline{\omega}}$, which could possibly be used
to cancel the term $\omega(X)$ in Equation (\ref{d1f1}). But
$\omega$ is locally exact, $\omega\!\!\mid_U=df_U$ ($U$: open
subset of $M$, $f_U\in\Ci(U)$). So it suffices to ensure that the
inner derivations associated to the functions $f_U$ glue together.
Lastly,
\[\overline{\omega}(D)\m_U=[D\m_U,f_U],\] $D\in\cD(M)$, and
$\Phi_2=\Phi_1\circ e^{-\kappa^{-1}\overline{\omega}}\in
Aut(\cD(M),[.,.])$ actually verifies
\begin{equation}\Phi_2(f+X)=\kappa f+\La\,div\,X+X.\label{Phi2}\end{equation} It is
interesting to note that the automorphism $e^{\overline{\omega}}$
is, for $\omega=df$ ($f\in\Ci(M)$), simply the inner automorphism
$e^{\overline{\omega}}:\cD(M)\ni D\raa e^f\cdot D\cdot
e^{-f}\in\cD(M)$.

An analogous extension $\overline{div}$ of the cocycle $div$ is at
least not canonical.

At this stage a new idea has to be injected. It is easily checked
that every quantum automorphism $\Phi$ of $\cD(M)$ induces a
classical automorphism $\tilde{\Phi}$ of $\cS(M)$,
$\tilde{\Phi}(\sigma(D))=\sigma(\Phi(D))$. Of course, the converse
is not accurate. Moreover, any automorphism of $\cS(M)$ restricts
to an automorphism of $\cS_0(M)\oplus\cS_1(M)\simeq\cD^1(M).$ So
it would be natural to try to benefit from this "algebra
hierarchy" and compute, having already obtained the automorphisms
of $\cD^1(M)$, those of $\cS(M).$ This approach however turns out
to be rather merely elegant than really necessary. We do not
employ it here.

The automorphism property shows that
\begin{equation}\Phi_2\m_{\cD^{i}(M)}=\kappa^{1-i}id+\psi_i,\label{structure}\end{equation} where $i\in\N$
and $\psi_i\in Hom_{\R}(\cD^{i}(M),\cD^{i-1}(M))$. This is
equivalent to saying that
\[\tilde{\Phi}_2\m_{\cS_i(M)}=\kappa^{1-i}id.\]
Such a classical automorphism ${\cal U}_{\kappa}:\cS_i(M)\ni P\raa
\kappa^{1-i}P\in\cS_i(M)$ really exists. Indeed, the degree
$deg:\cS_i(M)\ni P\raa (i-1)P\in\cS_i(M)$ is known to be a
derivation of $\cS(M)$ and, for $\kappa>0$, the automorphism of
$\cS(M)$ generated by $-log(\kappa)\cdot deg$ is precisely ${\cal
U}_{\kappa}$. Nevertheless, it is hard to imagine that this
classical automorphism ${\cal U}_{\kappa}$ is induced by a quantum
automorphism, since no graduation exists on the quantum level.
Therefore, the first guess is \[\kappa\stackrel{?}{=}1.\]

In order to validate or invalidate this supposition, we project
the automorphism property, combined with the newly discovered
structure (\ref{structure}) of $\Phi_2$, onto $\Ci(M).$ It is
worth stressing that all information is, here as well as above,
enclosed in the $0$-order terms.

Our exploitation of this projection is based on formal calculus.
In the main, this symbolism consists in the substitution of
monomials $\xi_1^{\ap^1}\ldots\xi_n^{\ap^n}$ in the components of
a linear form $\xi\in (\R^n)^*$ to the derivatives
$\p_{x^1}^{\ap^1}\ldots\p_{x^n}^{\ap^n}f$ of a function $f$.

This method is known in Mechanics as the \textit{normal ordering}
or \textit{canonical symbolization/quantization}. Its systematic
use in Differential Geometry is originated in papers by M. Flato
and A. Lichnerowicz \cite{FL} as well as M. De Wilde and P.
Lecomte \cite{DWLc2}, dealing with the Chevalley-Eilenberg
cohomology of the Lie algebra of vector fields associated with the
Lie derivative of differential forms. This polynomial modus
operandi matured during the last twenty years and developed into a
powerful computing technique, successfully applied in numerous
works (see e.g. \cite{LMT} or \cite{NP2}).

A by now standard application of the normal ordering leads to a
system of equations in specially the above constants $\La$ and
$\kappa$. When solving the system, we get two possibilities,
$(\La,\kappa)=(0,1)$ and $(\La,\kappa)=(1,-1)$. This outcome
surprisingly cancels the conjecture $\kappa=1.$\\

2.5.2.1 $(\La,\kappa)=(0,1)$\\

Equation (\ref{Phi2}), which gives $\Phi_2$ on $\cD^1(M)$,
suggests that $\Phi_2$ could coincide with $id$ on the whole
algebra $\cD(M).$

In fact our automorphism equations show that computations reduce
to the determination of some intertwining operators between the
$\cX(M)$-modules of $k$th-order linear differential operators
mapping differential $p$-forms into functions. These equivariant
operators have been obtained in \cite{NP} and \cite{BHMP}. They
allow to conclude that we actually have $\Phi_2=id$.

The paper \cite{NP}, following works by P. Lecomte, P. Mathonet,
and E. Tous\-set \cite{LMT}, H. Gargoubi and V. Ovsienko
\cite{GO}, P. Cohen, Yu. Manin, and D. Zagier \cite{CMZ}, C. Duval
and V. Ovsienko \cite{DO}, gives the classification of the
preceding modules. Additionally, it provides the complete
description of the above-mentioned intertwining operators, thus
answering a question by P. Lecomte whether some homotopy
operator---which locally coincides with the Koszul differential
\cite{PLc} and is equivariant if restricted to low-order
differential operators---is intertwining for all orders of
differentiation.

A small dimensional hypothesis in \cite{NP}, which was believed to
be inherent in the used canonical symbolization technique, was the
starting point of \cite{BHMP}. Here, the authors prove the
existence and the uniqueness of a projectively equivariant symbol
map (in the sense of P. Lecomte and V. Ovsienko \cite{PLVO}) for
the spaces of differential operators transforming $p$-forms into
functions, the underlying manifold being endowed with a flat
projective structure. The substitution of this equivariant symbol
to the previously used canonical symbol allowed to get rid of the
dimensional assumption and unexpected intertwining operators were
discovered in the few supplementary dimensions.\\

2.5.2.2 $(\La,\kappa)=(1,-1)$\\

Computations being rather technical, we confine ourselves to
mention that our quest for automorphisms, by means of the
canonical symbolization, leads to a symbol that reminds of the
opposite of the conjugation operator. Thus, this operator might be
an astonishing automorphism. Remember that for an oriented
manifold $M$ with volume $\Omega$, the conjugate $D^*\in\cD(M)$ of
a differential operator $D\in\cD(M)$ is defined by
\[\int_MD(f)\cdot g\mid\Omega\mid=\int_Mf\cdot D^*(g)\mid\Omega\mid,\]
for any compactly supported $f,g\in\Ci(M)$. Since
$(D\circ\DE)^*=\DE^*\circ D^*$, $D,\DE\in\cD(M)$, the operator
$\cC:=-*$ verifies $\cC(D\circ\DE)=-\cC(\DE)\circ\cC(D)$ and is
thus an automorphism of $\cD(M)$. Formal calculus allows to show
that this automorphism exists for any manifold (orientable or
not). At last, the computations can again be reduced to the
formerly described intertwining operators and we get $\Phi_2=\cC$.
Hence the conclusion:
\begin{theo} A linear map $\Phi:\cD(M)\raa\cD(M)$ is  an
automorphism of the Lie  algebra $\cD(M)$ of linear differential
operators on $C^\infty(M)$ if and only if it can be written in the
form \[\Phi=\varphi_*\circ\cC^a\circ e^{\overline{\omega}},
\]
where $\varphi$  is  a  diffeomorphism  of  $M$, $a=0,1$,
$\cC^0=id$ and $\cC^1=\cC$, and $\omega$  is  a  closed 1-form on
$M$. All the objects $\varphi,a,\omega$ are uniquely determined by
$\Phi$.
\end{theo}

\subsubsection{Automorphisms of $\mathbf{\cS(M)}$}

The study of the automorphisms of $\cS(M)$ is similar to the
preceding one regarding $\cD(M)$ and even simpler, in view of the
existence of the degree-automorphism ${\cal U}_{\kappa}$. We
obtain the following upshot.
\begin{theo} A linear map $\Phi:\cS(M)\raa \cS(M)$ is an automorphism
of the Lie algebra $\cS(M)$ of polynomial functions on $T^*M$ with
respect to the canonical symplectic bracket if and only if it can
be written in the form \[ \Phi(P)={\cal
U}_{\kappa}(P)\circ\varphi^*\circ Exp(\omega^v),
\]
where $\kappa$ is a non-zero  constant, ${\cal
U}_{\kappa}(P)=\kappa^{1-i}P$ for $P\in\cS_i(M)$, $\varphi^*$ is
the phase lift  of  a diffeomorphism $\varphi$  of  $M$,  and
$Exp(\omega^v)$ is the vertical symplectic diffeomorphism of
$T^*M$, which is nothing but the translation by a closed 1-form
$\omega$ on $M$. All  the objects $\kappa,\varphi,\omega$ are
uniquely determined by $\Phi$. \label{SM}\end{theo} It is
interesting to compare this result with those of \cite{AG} and
\cite{JG2}.

\bigskip
\textit{Acknowledgement.} The authors are grateful to Peter Michor for
providing helpful comments on topologies in spaces of differentiable maps.

\noindent Janusz GRABOWSKI\\Polish Academy of Sciences, Institute
of Mathematics\\\'Sniadeckich 8,
P.O.Box 21, 00-956 Warsaw, Poland\\Email: jagrab@impan.gov.pl\\\\\
\noindent Norbert PONCIN\\University of Luxembourg, Mathematics Laboratory\\avenue de la Fa\"{\i}encerie, 162 A\\
L-1511 Luxembourg City, Grand-Duchy of Luxembourg\\Email:
norbert.poncin@uni.lu

\begin{thebibliography}{Dillo 83}
\bibitem[Abe82]{A} Abe K, \textit{Pursell-Shanks type theorem for orbit spaces and
$G$-manifolds}, Publ. Res. Inst. Math. Sci., \textbf{18} (1982),
pp. 265-282
\bibitem[Ame75]{Am} Amemiya I, \textit{Lie algebra of vector fields and complex
structure}, J. Math. Soc. Japan, \textbf{27} (1975), pp. 545-549
\bibitem[AG90]{AG} Atkin C J, Grabowski J, \textit{Homomorphisms of the Lie
algebras associated with a symplectic manifold}, Compos. Math.,
\textbf{76} (1990), pp. 315-348
\bibitem[BHMP02]{BHMP} Boniver F, Hansoul S, Mathonet P, Poncin N,
\textit{Equivariant symbol calculus for differential operators
acting on forms}, Lett. Math. Phys., \textbf{62} (2002), pp.
219-232
\bibitem[CGM99]{CGM} Campillo A, Grabowski J, M\"uller G,
\textit{Derivation algebras of toric varieties}, Comp. Math.,
\textbf{116} (1999), pp. 119-132
\bibitem[CMZ97]{CMZ} Cohen P, Manin Yu, Zagier D, \textit{Automorphic pseudodifferential
operators, Algebraic Aspects of Integrable Systems}, Progr.
Nonlinear Differential Equations Appl., \textbf{26} (1997), pp.
17-47, Birkh\H{a}user Verlag
\bibitem[DWL81]{DWLc1} De Wilde M, Lecomte P, \textit{Some Characterizations of
Differential Operators on Vector Bundles}, E.B. Christoffel,
Butzer P, Feher F (ed), Brikh\H{a}user Verlag, Basel (1981), pp.
543-549
\bibitem[DWL83]{DWLc2} De Wilde M, Lecomte P, \textit{Cohomology of the
Lie Algebra of Smooth Vector Fields of a Manifold, associated to
the Lie Derivative of Smooth Forms}, J. Math. pures et appl.,
\textbf{62} (1983), pp. 197-214
\bibitem[DO97]{DO} Duval C, Ovsienko V, \textit{Space of second order
linear differential operators as a module over the Lie algebra of
vector fields}, Adv. in Math., \textbf{132}(2) (1997), pp. 316-333
\bibitem[FL80]{FL} Flato M, Lichnerowicz A, \textit{Cohomologie des repr\'{e}sentations definies par la derivation de
Lie et \`{a} valeurs dans les formes, de l'alg\`{e}bre de Lie des champs
de vecteurs d'une vari\'{e}t\'{e} diff\'{e}rentiable. Premiers espaces de
cohomologie. Applications.}, C. R. Acad. Sci., S\'{e}r. A \textbf{291}
(1980), pp. 331-335
\bibitem[GO96]{GO} Gargoubi H, Ovsienko V, \textit{Space of linear
differential operators on the real line as a module over the Lie
algebra of vector fields}, Internat. Math. Res. Notices,
\textbf{5} (1996), pp. 235-251
\bibitem[GG01]{GG} Grabowska K, Grabowski J, \textit{The Lie algebra of a
Lie algebroid}, in "Lie Algebroids and Related Topics in
Differential Geometry", J.~Kubarski et al. (Eds.), Banach Center
Publications, \textbf{54}, Warszawa 2001, pp. 43-50
\bibitem[GK37]{GK} Gel'fand I, Kolmogoroff A, \textit{On rings of continuous
functions on topological spaces}, C. R. (Dokl.) Acad. Sci. URSS,
\textbf{22} (1939), pp. 11-15
\bibitem[Gra78]{JG} Grabowski J, \textit{Isomorphisms and ideals of the Lie
algebras of vector fields}, Invent. math., \textbf{50} (1978), pp.
13-33
\bibitem[Gra93]{JG3} Grabowski J, \textit{Lie algebras of vector fields and generalized
foliations}, Publ. Matem., \textbf{37} (1993), pp 359-367
\bibitem[Gra00]{JG2} Grabowski J, \textit{Isomorphisms of Poisson and Jacobi
brackets}, in Poisson Geometry, Eds: J. Grabowski and P.
Urba\'nski, Banach Center Publications, \textbf{51}, Warszawa
2000, pp. 79-85
\bibitem[Gra03]{Gra03} Grabowski J, \textit{Isomorphisms of algebras of
smooth functions revisited}, Archiv Math. (to appear) (electronic
version at http://arXiv.org/abs/math.DG/0310295)
\bibitem[GP03]{GP} Grabowski J, Poncin N, {\it Automorphisms of quantum and classical Poisson algebras},
Comp. Math., London Math. Soc., 140, 2004, pp 511-527
\bibitem[Gt58]{Grau} Grauert H, \textit{On Levi's problem and the
embedding of real analytic manifolds}, Ann. Math., \textbf{68}
(1958), pp. 460-472
\bibitem[HM93]{HM} Hauser H, M\"{u}ller G, \textit{Affine varieties and Lie algebras of vector
fields}, Manusc. Math., \textbf{80} (1993), pp. 309-337
\bibitem[Kor74]{Ko} Koriyama A, \textit{On Lie algebras of vector fields
with invariant submanifolds}, Nagoya Math. J., \textbf{55} (1974),
pp. 91-110
\bibitem[KMO77]{KMO} Koriyama A, Maeda Y, Omori H, \textit{On Lie algebras
of vector fields}, Trans. Amer. Math. Soc., \textbf{226} (1977),
pp. 89-117
\bibitem[KM97]{KM} Kriegl A, Michor P W, \textit{The Convenient Setting of
Global Analysis}, Math. Surv. Monog., \textbf{53}, American
Mathematical Society 1997
\bibitem[Lec94]{PLc} Lecomte P, \textit{On some sequence of graded Lie algebras associated to
manifolds}, Ann. Glob. Anal. Geom., \textbf{12} (1994), pp. 183-192
\bibitem[LMT96]{LMT} Lecomte P, Mathonet P, Tousset E, \textit{Comparison
of some modules of the Lie algebra of vector fields}, Indag.
Math., \textbf{7(4)} (1996), pp. 461-471
\bibitem[LO99]{PLVO} Lecomte P, Ovsienko V, \textit{Projectively
equivariant symbol calculus}, Lett. Math. Phys., \textbf{49}
(1999), pp. 173-196
\bibitem[Mr\v c03]{Mrc} Mr\v cun J, \textit{On isomorphisms of algebras of
smooth functions}, electronic version at
http://arXiv.org/abs/math.DG/0309179
\bibitem[Omo76]{O} Omori H, \textit{Infinite dimensional Lie transformation
groups}, Lect. Notes in Math., \textbf{427} (1976), Springer
Verlag
\bibitem[Pon99]{NP2} Poncin N, \textit{Cohomologie de l'alg\`{e}bre de Lie des op\'{e}rateurs diff\'{e}rentiels sur une vari\'{e}t\'{e}, \`{a} coefficients
dans les fonctions}, C.R.A.S. Paris, \textbf{328} S\'{e}rie I
(1999), pp. 789-794
\bibitem[Pon02]{NP} Poncin N, {\it Equivariant Operators between some Modules of the Lie Algebra
of Vector Fields}, Comm. in Alg., 32 (7), 2004, pp 2559-2572
\bibitem[PS54]{PS} Shanks M E, Pursell L E, \textit{The Lie algebra of a smooth
manifold}, Proc. Amer. Math. Soc., \textbf{5} (1954), pp. 468-472
\bibitem[Sie96]{Si} Siebert T, \textit{Lie algebras of derivations and
affine differential geometry over fields of characteristic 0},
Mat. Ann., \textbf{305} (1996), pp. 271-286
\bibitem[Skr87]{S} Skryabin S M, \textit{The regular Lie rings of derivations of commutative
rings}, preprint WINITI 4403-W87 (1987)
\bibitem[Whi36]{Wh} Whitney H, \textit{Differentiable manifolds},
Ann. Math., \textbf{37} (1936), pp. 645-680
\end{thebibliography}
\end{document}